\documentclass[12pt]{article}
\usepackage{color}
\newtheorem{thm}{\indent{\sc Theorem}}
\newtheorem{lem}{\indent{\sc Lemma}}
\newtheorem{cor}{\indent{\sc Corollary}}

\newtheorem{rmk}{\indent{\sc Remark}}

\newcommand{\be}{\begin{eqnarray}}
\newcommand{\beo}{\begin{eqnarray*}}
\newcommand{\ee}{\end{eqnarray}}
\newcommand{\eeo}{\end{eqnarray*}}
\newcommand{\fa}{\frac{1}{\alpha}}
\newcommand{\fpa}{\frac{p}{\alpha}}

\newcommand{\n}{\nonumber}

\newcommand{\y}{y^{\frac{2}{\alpha}}}

\newcommand{\nn}{\nonumber}

\newcommand{\qed}{\hfill \rule{1.3mm}{3mm}}

\begin{document}

\title{{Process convergence of self normalized sums of i.i.d. random variables coming from domain of attraction of stable distributions}}
\author{ G K Basak \footnote{Stat-Math Unit, Indian Statistical Institute, Kolkata}    \& Arunangshu Biswas \footnote{Dept. of Statistics, Presidency College, Kolkata }}
\date{ }
\maketitle

\begin{abstract}
In this paper we show that the continuous version of the self normalised process 
$Y_{n,p}(t)= S_n(t)/V_{n,p}+(nt-[nt])X_{[nt]+1}/V_{n,p}$ where $S_n(t)=\sum_{i=1}^{[nt]} X_i$ and 
$V_{(n,p)}=(\sum_{i=1}^{n}|X_i|^p)^{\frac{1}{p}}$ and $X_i$ i.i.d. random variables belong to $DA(\alpha)$,
 has a non trivial distribution iff $p=\alpha=2$. The case
 for $2 > p > \alpha$ and $p \le \alpha < 2$  
 is systematically eliminated by showing that either of tightness or finite
  dimensional convergence to a non-degenerate limiting distribution does not hold. This work is an extension 
  of the work by Cs\"org\"o et al. who showed Donsker's theorem for 
     $Y_{n,2}(\cdot)$, i.e., for $p=2$, holds iff $\alpha =2$
  and identified the limiting
   process as standard Brownian motion in sup norm. 
\end{abstract}

\section{Introduction}
\label{sec.1}
\setcounter{equation}{0}
Limit theory plays a fundamental role in probability and statistics. Various forms of limit theorems, like
 the strong laws of large numbers, the central
limit thoerems, the law of iterarted logrithm and the laws of large deviations
are celebrated results in this field. However restrictive assumptons like the finiteness of moments upto a 
certain order or the existence of the moment generating function in
a neighbourhood of zero are necessary conditions for proving these
theorems. Also the choice of the normalising factor is the standard deviation,
which is typically unknown in many statistical applications. What is done instead is to estimate the unknown
 parameters by a sequence of random
variables ( the sample standard deviation like the Student's t statistic). The normalising factor is random in this case.
To see whether the above mentioned limit laws hold with random normalisation is a fruitful area of research that has yielded
many interesting results in the last two decades. For example, it has been shown in \cite{Kuelbs} that even
 under much less
assumptions an analogy of the law of iterated logarithm holds under randomised normalisation. The same thing
 can be shown in case of laws of large and moderate deviations, see \cite{Shao}.\\
The study of the asymptotics of the self normalised sums are also interesting. Logan et al, \cite{Logan} 
first showed the asymptotics of the self normalised sums where the variables belong to the domain of 
attraction of a stable distribution. In \cite{Gine}, it has been shown that limiting distribution of the 
self normalised sums converges to Normal if and only if the constituent random variables coming 
from the domain of attraction of a Normal stable distribution (henceforth denoted as $DAN$).
Hence they conclude the same for t-statistics.
Cs\"org\"o et al \cite {Csorgo} show a functional (process) convergence result 
in sup norm for suitably scaled products of the self normalised sums (with $L_2$ normalisation as in \cite{Gine}).
They also show the result holds if and only if constituent random variables come
  from $DAN$. 
Basak et al 
  \cite {Basak} showed the convergence of a suitably scaled process to an Ornstein Uhlenbeck process. There
   also the constituent variables come from $DAN$. The aim of this paper is to show that the only case when the asymptotic 
   distribution of the self normalised process is non trivial is when the norming index $p$ is exactly equal to the index of stability $\alpha$  (for definition see Section \ref{defn}).\\
This paper is organised as follows. Section \ref{defn} contains definitions and a preliminary result that is used throughout.
Section \ref{main} contains the main result of this paper. Section \ref{application} contains  various application applied to functionals of the self normlaised process
  as corollary to the main theorem. Section \ref{finite} and Section \ref{tight} shows
 convergence of finite dimensional distribution of the self normlaised process and tightness result 
 respectively for various choices of $p$ and $\alpha$. This two sections together show what should be the 
 relation between $p$ and $\alpha$ for which the resulting asymptotic distribution is non trivial. Secton \ref{conclusion} concludes the paper with a few possible research directions.

\section {Definition and preliminaries}
\label{defn}
\setcounter{equation}{0}

Let $\{X_i\}$ be a sequence of i.i.d. random variables.
We intend to study the convergence of the process determined at time $t$ by \\
\begin{eqnarray}
\label{defn eq}
Y_{n,p}(t)=\frac{S_n(t)}{(V_{n,p})} + (nt-[nt]) X_{[nt]+1}/V_{n,p}& & \mbox{$0<t<1$}   
 \; \; \; \; \mbox{$p > 0$}
\end{eqnarray}
where the process $S_n(.)$ and $V_{n,p}(.)$ is defined as\\

$
S_n(t)=
\sum_{i=1}^{[nt]} X_i 
$
and
$V_{n,p}= 
(\sum_{i=1}^{n} |X_i|^{p})^{1/p} 
$
where $X_i$'s belong to the domain of attraction of a $\alpha$-stable family denoted by $DA(\alpha)$ and
$[x]$ is the largest integer less than or equal to $x$. We
prove process convergence by showing finite dimensional convergence and tightness.

\bigskip

Here we prove a lemma first for the benefit of the reader, as we would use it often (Feller, Vol. 2, \cite{Feller}) :\\
\begin{lem}
\label{lem.prelim}
If $X \in DA(\alpha)$ then $Y= sgn(X)|X|^{\frac{\alpha}{2}} \in DAN$.
\end{lem} 

\textbf{Proof.} To prove the lemma we need the following characterisation:
\begin{center}
$Y \in DAN$ iff $\lim_{y \rightarrow \infty}\frac{y^2 P(|Y|>y)}{E(Y^2 I(|Y|<y))} = 0$ ,
\end{center}
see \cite{Csorgo}\\
We show that the random variable Y satisfies the necessary and sufficient condition.\\
Now, 

\begin{eqnarray*}
\lim_{y \rightarrow \infty}y^2P(|Y|>y)&=&\lim_{y \rightarrow \infty}y^2 P(|X|^{\frac{\alpha}{2}}>y)\\
&=& \lim_{y \rightarrow \infty} y^2 P(|X|>y ^{\frac{2}{\alpha}})\\
&=& \lim_{y \rightarrow \infty} y^2 (y^{\frac{2}{\alpha}})^{-\alpha}\;\;  \\
&&  \mbox{since the tail of DA($\alpha$) is Paretian, i.e, $P(|X|>x) = O(x^{-\alpha})$}\\
&=& O(1) .
\end{eqnarray*}

And,
\begin{eqnarray*}
E(Y^2 I(|Y|<y))&=& E(|X|^{\alpha} I(|X|^{\frac{\alpha}{2}}\le y))\\
&=& E(|X|^{\alpha} I (|X|\le y^{\frac{2}{\alpha}}))\\
&=& \int_{0}^{y^{\frac{2}{\alpha}}} z^{\alpha} d F_{|X|}(z)\\
&=& \int_{0}^{\y} (\int_{0}^{z} \alpha t ^{\alpha -1} dt ) dF_{|X|}(z) .\\
\end{eqnarray*}
{Applying Fubini's theorem and interchanging order of integration we get }\\

\begin{eqnarray*}
E(Y^2 I(|Y|<y))
&=& \int_{0}^{\y} \alpha \int_{t}^{\y} dF_{|X|}(z) t^{\alpha -1} dt\\
&=& \alpha \int_{0}^{\y} P(t<|X| \le \y) t^{\alpha -1}dt\\
&=& \alpha \int_{0}^{M} P(t<|X| \le \y) t^{\alpha -1}dt\\
&+& \alpha \int_{M}^{\y} P(t<|X| \le \y) t^{\alpha -1}dt .\\
\end{eqnarray*}
Now the first integral is nonnegative and less than $M^\alpha$ and
 for the limit of the second integral as $y \rightarrow \infty $, 
 use Monotone Convergence Theorem, to get  
\begin{eqnarray*}
& & \alpha \int_{M}^{\y} P(t<|X|< \y) t^{\alpha -1}dt\\
&=&\lim_{y\rightarrow \infty} \alpha \int_{M}^{\infty} I_{(M, \y]}(t) \ P(t<|X|\le \y) t^{\alpha - 1} dt\\
&=& \alpha \int_{M}^{\infty}  P(|X| > t) t^{\alpha - 1} dt\\
&=& \infty , \ \ \mbox{ as } P(|X| > t) = O(t^{-\alpha})  .
\end{eqnarray*}

Hence, $\lim_{y \rightarrow \infty}\frac{y^2 P(|Y|>y)}{E(Y^2 I(|Y|<y))} = 0 \Leftrightarrow Y \sim DAN.$
\qed

\bigskip
We also quote a theorem due to \cite{Csorgo}
\begin{thm}
\label{thm.1}
The following statements are equivalent: 

\begin{enumerate}
	\item
$EX = 0$ and X is in the domain of attraction of the normal law. 
\item
$S_{[nt_0]}/ V_{n,2} \to N(0, t_0)$ for $t_0 \in (0, 1]$. 
\item
$S_{[nt]}/V_n \to W(t)$ on $(D[0, 1], \rho)$, where $\rho$ is the sup-norm metric for functions in $D[0, 1]$, and $\{W(t), 0 < t < 1\}$ is a standard Wiener process. 
\item
On an appropriate probability space for $X, X_1, X_2, \ldots$ we can construct a standard Wiener process $\{W(t), 0 < t < \infty\}$ such that 
\be
\sup_{0\le t \le 1} |S_{[nt]}/ V_{n,2} -W(nt)/\sqrt{n}| &=& o_p(1).   \label{supnorm}
\ee
\end{enumerate}
\end{thm}

\section{Main result}
\label{main}
\setcounter{equation}{0}

Let $X_i$ be i.i.d. symmetric observations from the domain of attraction of a $\alpha$-Stable distribution and 
$\{Y_{n,p}(\cdot)\}$ as defined in (\ref{defn eq}). Then we have the following theorem:
\begin{thm}
$Y_{n,p}(t)$ converges weakly to Brownian motion in $C[0,1]$, if and only if $p = \alpha = 2$.
\end{thm}
\textbf{Proof:} 
In Section \ref{finite} we show that for $0 < p<\alpha \le 2$ and $0 < p=\alpha<2$ the
finite dimensional distributions converge in probability to a degenerate distribution at zero. 
A non trivial limiting distribution exists if $p > \alpha$ and $p = \alpha = 2$. In sections \ref{tight} we show that the 
sequence $\{S_n/V_{n,p}\}$ of self normalised sums is tight iff $0< p \le \alpha \le 2$. The only 
case where we have both tightness and finite dimensional convergence is $p=\alpha=2$. The limiting
 distribution of the  sequence for this choice of $p$ and $\alpha$ was identified by \cite{Gine} as Normal. 
 Applying Prohorov's Theorem we have the distributional convergence to the Wiener process. The convergence 
 in the sup norm metric follows directly from (\ref{supnorm}) of the above theorem by \cite{Csorgo}.\\
 
 {In Cs\"org\"o (\cite{Csorgo}) they are interested in the process $S_{[nt]}/V_{n,p}$ which is in $D([0,1])$ 
 and we are interested in the process $Y_{n,p}(t)$ which is in $C([0,1])$. However from the definition of
  $Y_{n,p}(t),$
 \beo
  |Y_{n,p}(t)- S_{[nt]}/V_{n,p}| &=& |(nt-[nt])X_{[nt]+1}|/V_{n,2} \le |X_{[nt]+1}|/V_{n,p}.
  \eeo
  If \ $p=\alpha=2$ \ then, by Darling (\cite{Darling}), we have that 
  \ $\max_{1\le i \le n} |X_i|/(\sum_{i=1}^n X_i^2)^\frac{1}{2} \stackrel{P}\longrightarrow 0.$ 
  \ So \ $|Y_{n,p}(t)- S_{[nt]}/V_{n,p}| \stackrel{P}\rightarrow 0$. Therefore $Y_{n,p}(t)$ takes the same
  limiting distribution of $S_{[nt]}/V_{n,p}$ which is Normal.}

\section{Application}
\label{application} 
\setcounter{equation}{0}

Here we present a few applications of the main theorem. These follows from  the original extension by Erdos and Kac \cite {erdos} to the corresponding self normalised functionals. Let $W(t)$ be the standard Brownian motion. Define the following quantities\\
$
\begin{array}{ll}
G_1(x)= P\Bigl(\sup_{0<t<1} W(t) < x \Bigr),   &   \;\;\; G_2(x)= P\Bigl(\sup_{0<t<1} |W(t)| < x \Bigr) \\
G_3(x)= P\Bigl( 	\int_{0}^{1} W^2(t) < x \Bigr),  &  \;\;\; G_4(x)= P\Bigl( 	\int_{0}^{1} |W(t)| < x \Bigr)
\end{array}
$

\begin{cor}
The following weak convergence holds iff $p=\alpha=2, \forall x>0 $.
\begin{enumerate}
	\item $P(\max_{1 \le k \le n} S_k/ V_{n,p}<x) \rightarrow G_1(x)$
	\item $P(\max_{1 \le k \le n}) |S_k|/V_{n,p} \rightarrow G_2(x)$
	\item $P(\frac{1}{n}\sum_{1 \le k \le n} (S_k/V_{n,p})^2) \rightarrow G_3(x)$
	 \item $P(\frac{1}{n} \sum_{1 \le k \le n} |S_k/V_{n,p}|) \rightarrow G_4(x)$
\end{enumerate}
\end{cor}

\section{Convergence of Finite Dimensional Distributions}
\label{finite}
\setcounter{equation}{0}

To get the process convergence we first need to examine the convergence of finite-dimensional
distributions, i.e.,
for $0 < t_1<t_2 <\ldots <t_k, k \ge 1$ we want to examine the convergence of the random vector
 $(Y_{n,p}(t_1), Y_{n,p}(t_2), \ldots, Y_{n,p}(t_k))$ as $n \rightarrow \infty .$ 
We will do this for $p < \alpha$, $p = \alpha$ and $p > \alpha$ separately.
\subsection {Case 1: $p < \alpha$}
Since $X_i \in DA(\alpha)$, 
by SLLN, \ $V_{n,p}/n^{1/p}$ converges to a positive constant, say, $k(\alpha, p)$.

Now, for $X_i \in DA(\alpha)$, \ $S_n/(n^{1/\alpha} h(n))$ converges in distribution to a $S(\alpha)$ random
variable,
where $h$ is a slowly varying function of $n$. 
Since $p < \alpha$, \ $S_n/n^{1/p} = n^{(1/\alpha) - (1/p)} S_n/n^{1/\alpha} \to 0$, in probability, as $n \to \infty$.
Thus, $S_n/V_{n,p} = \frac{S_n/n^{1/p}}{V_{n,p}/n^{1/p}} \to 0 , \ \mbox{ in probability, as } \  
n \to \infty .$
Therefore, the joint distribution would converge to a degenerate one, in this case.

\subsection {Case 2: $p=\alpha$}

Here we assume that 
$X_i$ is symmetric and belongs to $DA(\alpha)$.

\begin{lem}
For $V_{n,\alpha}$ defined as in Section \ref{defn} 
$V_{n,\alpha}\ge V_{n, 1} \ge V_{n,\beta} \ge V_{n,2}$ \ if \ $\alpha\le 1 \le \beta \le 2$.  
\end{lem}

\textbf{Proof.}
We use the inequality for $a>0, b>0$, and $\alpha \le 1, \ 2 \ge \beta \ge 1$,
\beo
a^\alpha+b^\alpha \ge (a+b)^\alpha  & \mbox{ and } &  (a+b)^\beta \ge a^\beta+b^\beta  \\
& \Rightarrow & (a^\alpha+b^\alpha)^{\frac{\beta}{\alpha}} \ge (a+b)^\beta 
\ge (a^\beta+b^\beta)\\
& \Rightarrow & (a^\alpha+b^\alpha)^{\frac{1}{\alpha}} \ge (a^\beta+b^\beta)^{\frac{1}{\beta}} .
\eeo
Now, take $\alpha = 1$ and $\beta \ge 1$ and then  $\alpha \le 1$ and $\beta = 1$ to get,
$$(a^\alpha+b^\alpha)^{\frac{1}{\alpha}} \ge (a+b) \ge (a^\beta+b^\beta)^{\frac{1}{\beta}} .$$  
Also, for $1 \le \beta \le 2$, it follows that $1 \le 2/\beta \le 2$. Hence,
$$  (a^\beta+b^\beta)^{2/\beta} \ge (a^{\beta (2/\beta)}+b^{\beta (2/\beta)}) = (a^2 + b^2) 
\Rightarrow (a^\beta+b^\beta)^{\frac{1}{\beta}} \ge (a^2+b^2)^{\frac{1}{2}} .
$$
The case for $n$ positive numbers can be shown in the same manner.
Thus, combining above, 
get
\beo
V_{n,2} \le V_{n,\beta} \le V_{n,1} \le V_{n,\alpha} .
\eeo
\qed

\bigskip
Now, we show that self-normalized sum for $p=\alpha$ converges to degenerate distribution as well.
\begin{thm}
\label{thm.3}
 If $p=\alpha \le 1,\; \lim_{n \rightarrow \infty} \mbox{Var} (\frac{S_n}{V_{n, p}})=0 .$ 
\end{thm}

\textbf{Proof.} 
 Note that,
\be
\label{eq.sym}
E\Bigl(\frac{\sum X_i}{V_{n, \alpha}}\Bigr)^2 &=& \sum E\Bigl(\frac{ X_i^2}{V_{n, \alpha}^{2}}\Bigr)\n 
+ \sum_{(i,j): i \neq j} E\Bigl(\frac{X_i X_j}{V_{n, \alpha}^{2}}\Bigr)\n \\
&=& \sum_i E\Bigl(\frac{X_i^2}{V_{n, \alpha}^{2}}\Bigr) + \sum_i E\Bigl(\sum _{j\neq i} X_i 
E\Bigl(\frac{X_j}{V_{n, \alpha}^{2}} \ | \ X_i, i \neq j \Bigr)\Bigr)\n \\
&=& \sum_i E \Bigl(\frac{X_i^2}{V_{n, \alpha}^{2}}\Bigr) , 
\ee
the second term vanishes since 
$$\frac{X_j}{V_{n, \alpha}^{2}} = - \frac{X_j}{V_{n, \alpha}^{2}}\;\; \mbox{ in distribution. } $$

We now
use the fact that \
$V_{n,\alpha} \ge V_{n,2}$ implies that $(\sum_{i=1}^n X_i^2/V_{n,\alpha})^2 \le 
(\sum_{i=1}^n X_i^2/V_{n,2})^2 = 1$.
Hence
if we could show that $(\sum_{i=1}^n X_i^2/V_{n,\alpha})^2 \to 0$ in probability, then
by Dominated convergence theorem (DCT) we have the result.
Observe that, for $\alpha \le 1$,
 if $X_i \sim DA(\alpha)$ \ 
 then $Y_i= sgn(X_i)|X_i|^\frac{\alpha}{2} \sim DAN$. 
 From \cite {Gine} we have that 
$E(\frac{Y_i^4}{(\sum |Y_i|^2)^2}) = E(\frac{X_i^{2\alpha}}{(\sum |X_i|^{\alpha})^2})= o(\frac{1}{n})$. 
 \ Thus,
\beo
E(\frac{\sum_{i=1}^n X_i^2 }{(V_{n, \alpha})^2})^{\alpha} 
& \le &
E(\frac{\sum_{i=1}^n |X_i|^{2 \alpha}}{(\sum_{i=1}^n |X_i|^{\alpha})^2}) \\
& = & o(1) \ \to  0 \ \ \mbox{ as } \ \ n \to 
\infty .
\eeo
Hence, $(\sum_{i=1}^n X_i^2/V_{n,\alpha})^2 \to 0$ in probability, as it goes to zero in $\alpha$-th mean.
Therefore, by DCT we conclude the proof.
\qed

\bigskip
We now proceed to prove the result for  $p=\alpha > 1$. 
\begin{lem}
\label{lem.convprob}
 If $X \sim DA(\alpha) $, then $\frac{\sum X_i^2}{V_{n,\alpha}^{2}} \stackrel{P} {\rightarrow} 0$.
\end{lem}

\textbf{Proof.} Observe that,
if $X_i\sim DAN$ then by \cite{Obrien} \ $\max_{1\le i\le n} \frac{|X_i|}{V_{n,2}} \stackrel{P}{\rightarrow} 
 0$. Because, if $X_i \sim DA(\alpha)$ then $Y_i= sgn(X_i)|X_i|^{\frac{\alpha}{2}} \sim DAN$ by the 
 lemma \ref{lem.prelim}. Therefore
\be
 \label{eq.obrien}
 \max_{1 \le i \le n} \frac{|Y_i|}{(\sum Y_i^2)^{\frac{1}{2}}} \stackrel{P}{\rightarrow} 0 
 & \Leftrightarrow & \ \max_{1 \le i \le n} \frac{|X_i|^{\frac{\alpha}{2}}}{(\sum |X_i|^\alpha)^{\frac{1}{2}}} \stackrel{P}{\rightarrow} 0 \nn\\
 \Leftrightarrow \ \max_{1 \le i \le n} \frac{|X_i|}{(\sum {|X_i|^\alpha})^{\frac{1}{\alpha}}} \stackrel{P}{\rightarrow} 0 & \Leftrightarrow &
\max_{1 \le i \le n} \frac{ |X_i|^2} {V_{n,\alpha}^{2}} \stackrel{P} {\rightarrow} 0 .
\ee
Again, from \cite{Feller,Darling}, since $X_i \sim DA(\alpha)$, one gets $|X_i|^2 \sim DA(\alpha/2)$. 
Define $Y_n^*= \max_{1\le i \le n} X_i^2$.
Hence, for $\epsilon, \eta> 0$ choose $\delta=\frac{\epsilon}{K_{\eta}}$ where $K_{\eta}$ is chosen so that
 $P(\sum X_i^2/Y_n^*> K_{\eta})< \eta/2$. (This is possible since by \cite{Darling} \ $Y_n^*/\sum X_i^2$ has 
 a limiting distribution and hence tight.)
\beo
P(\frac{\sum X_i^2}{V_{n,\alpha}^{2}} > \epsilon)&\le& P(\frac{\sum X_i^2}{V_{n,\alpha}^{2}} > \epsilon, \ \frac{Y_n^*}{V_{n, \alpha}^{2}} > \delta) 
\ + \ P(\frac{\sum X_i^2}{V_{n,\alpha}^{2}} > \epsilon, \ \frac{Y_n^*}{V_{n, \alpha}^{2}} \le \delta) \\
&\le& P(\frac{Y_n^*}{V_{n, \alpha}^{2}} > \delta) 
\ + \ P(\frac{\sum X_i^2}{X_n^*} \frac{Y_n^*}{V_{n, \alpha}^{2}} > \epsilon, \ \frac{Y_n^*}{V_{n, \alpha}^{2}} \le \delta) \\ 
&\le& P(\frac{Y_n^*}{V_{n,\alpha}^2} >  \delta) \ + \ P(\frac{\sum X_i^2}{Y_n^*} > \frac{\epsilon}{\delta}) .
\eeo

Choose $n_0$ sufficiently large so that the first probability is less than $\eta/2$. By the choice of 
$\delta$ we have the second probability less than $\eta/2$. Which implies that,

\beo
P(\sum X_i^2/V_{n,\alpha}^2> \epsilon) &<& \eta \;\;  \mbox{for $n \ge n_0$} .
\eeo

Hence the lemma is proved.
\qed 

\bigskip
\begin{thm}
\label{thm.4}
Let \ $1 < p=\alpha < 2$, \ and \ $X_i$s are symmetric and $X_i \sim DA(\alpha)$.
Then $\lim_{n \rightarrow \infty} \mbox{Var} (\frac{S_n}{V_{n, p}})=0 .$ 
\end{thm}
Note that $\mbox{Var} (\frac{S_n}{V_{n, p}})= E(\frac{S_n}{V_{n, p}})^2$ by symmetry of $X_i$
and also $E(\frac{S_n}{V_{n, p}})^2 = E(\frac{\sum X_i^2}{V_{n, p}^2})$ 
by (\ref{eq.sym}) in the proof of Theorem \ref{thm.3}. 
\beo
V_{n,\alpha} &\ge& V_{n,2} \;\; \ \ \mbox{for $0 < \alpha \le 2$}\\
\Rightarrow \frac{\sum X_i^2}{(\sum |X_i|^\alpha)^{\frac{2}{\alpha}}}&\le& \frac{\sum X_i^2}{\sum X_i^2} = 1 .\\
\eeo
Hence, by lemma \ref{lem.convprob} and applying bounded convergence theorem,
\beo
\lim_{n \rightarrow \infty}E(\frac{\sum X_i^2}{(\sum |X_i|^\alpha)^ {\
\frac{2}{\alpha}}}) = 0 .
\eeo
This proves the Theorem.

\begin{rmk}
For $X_i \sim DA(\alpha)$ symmetric, we, in fact, showed in theorems \ref{thm.3} and \ref{thm.4} that
$(S_n/V_{n,p}) \to 0$ in probability, for $0 < p = \alpha < 2$. Using same technique, it is immediate
that for any fixed \ $0 \le t \le 1$, \ $(S_{[nt]}/V_{n,p}) \to 0$, in probability, for \ $0 < p = \alpha < 2$ as well. 
The result for $k$ dimension can be obtained from the above result. 
 Note that the joint distribution of $(\frac{S_{[nt_1]}}{V_{n,p}}, \frac{S_{[nt_2]}}{V_{n,p}},\ldots,\frac{S_{[nt_k]}}{V_{n,p}} )$ can be 
 obtained from the joint distribution of $(\frac{S_{[nt_1]}}{V_{n,p}}, \frac{S_{[nt_2]}-S_{[nt_1]}}{V_{n,p}},\frac{S_{[nt_3]}-S_{[nt_2]}}{V_{n,p}}, \ldots,\frac{S_{[nt_k]}-S_{[nt_{k-1}]}}{V_{n,p}})$ by a linear transformation. We next show that the joint distribution of the latter converges to zero. 
  Write $S_1=\frac{S_{[nt_1]}}{V_{n,p}}, S_2=\frac{S_{[nt_2]}-S_{[nt_1]}}{V_{n,p}}$ and $S_k=\frac{S_{[nt_k]}-S_{[nt_{k-1}]}}{V_{n,p}}$. Now consider the 
  varinace of any linear combination of them  $V(a_1 S_1 + a_2 S_2 + \ldots + a_k S_k)$ where $a_i's$ are any arbitrary constants. 
  Due to independence the cross product term vanishes and by Theorem \ref{thm.4} the variances are zero which implies that any linear combination tends in probability to zero. Therefore 
  $\phi_{S_1,S_2,\ldots,S_k}(a_1,a_2,\ldots,a_k) \rightarrow 1$, where $\phi_{S_1, S_2, \ldots, S_k}$ is the characteristic function. Applying 
  continuity theorem we therefore have that the limiting joint distribution of $(S_1,S_2, \ldots, S_k)$ and hence $(S_{[nt_1]}/V_{n,p},S_{[nt_2]}/V_{n,p},\ldots,S_{[nt_k]}/V_{n,p})$is degenerate at 0. 
  \end{rmk}
  
  \subsection {Case 3: $p>\alpha$}
  \label{pgealpha}
We show that a limiting distribution exits in this case by finding the joint characteristic function of
 $(\frac{S_{m_1}}{V_{n,p}},\frac{S_{m_2}}{V_{n,p}},\ldots,\frac{S_{m_k}}{V_{n,p}})$ where
  $1 \le m_1 \le m_2 \le \ldots m_k \le n$. This is equivalent to finding the characteristic function of
   $(S_{m_1}/n^{\fa},(S_{m_2}-S_{m_1})/n^{\fa}, \ldots, (S_{m_k}-S_{m_{k-1}})/n^{\fa},V_{n,p}^p/n^{\fpa})$ by virtue of a transformation. 
   The characteristic function of the latter is
\beo \hspace{-20pt}
&& E(exp(i\frac{u_1}{n^{\fa}} S_{m_1})+ i \frac{u_2}{{n^{\fa}}}(S_{m_2}-S_{m_1})+ \ldots \\
&& + i \frac{u_k}{{n^{\fa}}}(S_{m_k}-S_{m_{k-1}}) + i \frac{l}{{n^{\frac{p}{\alpha}}}} V_{n,p}) \\
&=& E(exp(i\frac{u_1}{{n^{\fa}}} S_{m_1})+ \frac{t_2}{{n^{\fa}}} (S_{m_2}-S_{m_1})+ \ldots \\
&& + \frac{u_k}{{n^{\fa}}}(S_{m_k}-S_{m_{k-1}})+\frac{l}{{n^{\frac{p}{\alpha}}}} (V_{n,p}-V_{m_k,p}+ V_{m_k,p}-V_{m_{k-1},p} \\
&& + \ldots V_{m_2,p} - V_{m_1,p} + V_{m_1,p})))\\
&=& E(exp(i \{\frac{u_1}{{n^{\fa}}} S_{m_1}+ \frac{l}{{n^{\fpa}}} V_{m_1}\}\\
&& + \{ \frac{u_2}{{n^{\fa}}} (S_{m_2} - S_{m_1}) + \frac{l}{{n^{\fpa}}} (V_{m_2,p}-V_{m_1,p}) \}\\
&& + \ldots + \frac{l}{{n^{\fpa}}} (V_{n,p}^p - V_{{m_k},p}^p)))
\eeo 
Due to independence and identical distribution of $X's$ we have 
\beo
E(exp(i\frac{u_1}{{n^{\fa}}}S_{m_1}+ i\frac{l}{{n^{\fpa}}}V_{m_1,p})) &=& E(exp(iu_1\frac{X}{n^{\fa}}+ il \frac{|X|^p}{n^{\frac{p}{\alpha}}}))^{m_1}
\eeo
and
\beo
E(exp(i\frac{u_k}{{n^{\fa}}}(S_{m_k}- S_{m_{k-1}})+ i\frac{l}{{n^{\fpa}}}(V_{m_k,p}-V_{m_{k-1},p}))) &=& E(exp(i u_k \frac{X}{n^{\fa}}+ il \frac{|X|^p}{n^{\frac{p}{\alpha}}}))^{m_k-m_{k-1}}
\eeo
Now,
\beo
&& \{E(exp(^{iu\frac{X}{m_1^\fa}(\frac{m_1}{n})^\fa+iw{(\frac{|X|}{m_1^\fa})^p}(\frac{m_1}{n})^\frac{p}{\alpha}}))\}^{m_1} \\
&=& (\int exp({iu\frac{x}{m_1^\fa}(\frac{m_1}{n})^\fa+iw{(\frac{|x|}{m_1^\fa})^p}(\frac{m_1}{n})^\frac{p}{\alpha}} g(x) dx)^{m_1}\n \\
&=& (1+\int (exp({iu\frac{x}{m_1^\fa}(\frac{m_1}{n})^\fa +iw{(\frac{|x|}{m_1^\fa})^p}(\frac{m_1}{n})^\frac{p}{\alpha}}-1)) g(x) dx)^{m_1}\n \\
&=& (1+\frac{1}{m_1}\int (exp({iu y (\frac{m_1}{n})^\fa +iw |y|^p (\frac{m_1}{n})^\frac{p}{\alpha}})-1) g(m_1^\fa y) \n \\
&\times& {(m_1^\fa y)}^{\alpha+1}) \frac{dy}{y^{\alpha+1}})^{m_1} . \n\\
\eeo
Since $(exp({iu\frac{x}{m_1^\fa}(\frac{m_1}{n})^\fa +iw{(\frac{|x|}{m_1^\fa})^p}(\frac{m_1}{n})^\frac{p}{\alpha}})-1)$ is bounded by 2 and $m^\fa g(m^\fa y)$ is integrable we can apply DCT to get 
\beo
\lim_{m_1\rightarrow \infty}E(exp({iu\frac{X}{m_1^\fa}(\frac{m_1}{n})^\fa+iw{(\frac{|X|}{m_1^\fa})^p}(\frac{m_1}{n})^\frac{p}{\alpha}}))\}^{m_1}&=&\lim(1+\frac{c_{m_1,n}}{m_1}(u,w))^{m_1}\\
&=& exp({\lim c_{m_1,n}(u,w)}) ,
\eeo
where
\beo
c_{m_1,n}(u,w) &=& \int(exp({iuy(\frac{m_1}{n})^\fa+iw|y|^p(\frac{m_1}{n})^\fpa}) -1) g(m^\fa y) m^{\fa+1}dy.
\eeo
and  
\beo
\lim_{m_1,n}c_{m_1,n}(u,w) &=& \int (exp({iuy(t_1)^\fa+iw|y|^p(t_2)^\fpa}) -1)\frac{K(y)}{y^{\alpha+1}}dy ,\\
\eeo
where $K(y)=\lim_{m\rightarrow \infty}(m^\fa y)^{\alpha+1} g(m^\fa y)$ which is \\
$
K(y)=\left\{
\begin{array}{ll}
r & \mbox{if $y>0$}\\
s & \mbox{if $y<0$} 
\end{array}
\right.
$
\\
The same thing can be done for $E(exp(iu_k\frac{X}{n^{\fa}}+ il \frac{|X|^p}{n^{\frac{p}{\alpha}}}))^{m_k-m_{k-1}}$ and let us call it $c_{m_{k-1},m_k,n}(u_k,l)$. The limiting chracteristic function of  $(S_{m_1}/n^{\fa},(S_{m_2}-S_{m_1})/n^{\fa}, \ldots, (S_{m_k}-S_{m_{k-1}})/n^{\fa},V_{n,p}^p/n^{\fpa})$ is therefore (the limits are meant in such a way that $\frac{m_i}{n}$ tends to a constant as $m_i,n \rightarrow \infty$).
\beo
&& \lim_{m_1,m_2,\ldots,m_k,n \rightarrow \infty} E(exp(i(\frac{t_1}{n^{\fa}} S_{m_1})+ \frac{t_2}{{n^{\fa}}} (S_{m_2}-S_{m_1})+ \ldots \\
&+& \frac{t_k}{{n^{\fa}}} (S_{m_k}-S_{m_{k-1}}) +\frac{l}{{n^{\frac{p}{\alpha}}}} V_{n,p})) \\
 &=& \lim_{m_1,n \rightarrow \infty} c_{m_1,n} \lim_{m_2,m_1,n \rightarrow \infty}c_{m_1,m_2,n} \ldots\\
 & & \times \ \lim_{m_{k-1},m_k,n \rightarrow \infty} c_{m_{k-1},m_k}\\
 & & \times \ \lim_{m_k,n}E(exp(e^{il(\frac{|X|}{n^{\fa}})^{p}}))^{n-m_k} ,
\eeo
since $X\sim DA(\alpha)$ the last limit exists. Hence the finite-dimensional distribution converges
to some (non-trivial) non-degenerate distribution for $p > \alpha$. 

\begin{rmk}
For $p=\alpha=2$ the finite dimensional distribution of can be obtained by using the fact (see \cite{Gine})
 that $\frac{S_n}{\sqrt{n l(n)}} \stackrel{D}{\rightarrow} N(0,1)$ and $\frac{1}{nl(n)}V_{n,2}^2 \stackrel{P}{\rightarrow} 1$ 
 for a slowly varying function $l(\cdot)$. Applying the same argument as above we see that the distribution
  of $(Y_{n,p}(t_1), Y_{n,p}(t_2), \ldots, Y_{n,p}(t_k))$ can be obtained from the distribution of \\
  $(\frac{S_{[nt_1]}}{V_{n,p}}, \frac{S_{[nt_2]}-S_{[nt_1]}}{V_{n,p}}, \ldots, \frac{S_{[nt_k]}-S_{[nt_{k-1}]}}{V_{n,p}})$ 
  by a linear transformation. Now the components in the latter are uncorrelated and hence 
\beo
\frac{S_{[nt_1]}}{V_{n,p}}= \frac{\sqrt{[nt_1]l([nt_1])}}{\sqrt{nl(n)}} \frac{\frac{1}{\sqrt{[nt_1]l([nt_1])}}S_{[nt_1]}}{\frac{1}{\sqrt{n l(n)}}V_{n,p}}&\stackrel{D}{\rightarrow}& {\sqrt{t_1}}N(0,1), 
\eeo
( by using Slutsky's Theorem and the fact that $l(\cdot)$ is a sowly varying function).
Similar thing can be done for $\frac{S_{[nt_2]}-S_{[nt_1]}}{V_{n,p}}$ and the limiting distribution in that case will be ${\sqrt{t_2-t_1}} N(0,1)$. If $t_i < t_j$ then $Cov(\frac{S_{[nt_i}}{V_{n,p}},\frac{S_{[nt_j]}}{V_{n,p}}) = Cov(\frac{S_{[nt_i]}}{V_{n,p}},\frac{S_{[nt_j]}-S_{[nt_i]}+S_{[nt_i]}}{V_{n,p}})=V(\frac{S_{[nt_i]}}{V_{n,p}})= t_i = \min (t_i,t_j) $. Since the Jacobian of the transformation is one the finite dimensional distribution of $(Y_{n,p}(t_1),Y_{n,p}(t_2),\ldots, Y_{n,p}(t_k))$ is a multivariate normal distribution with dispersion matrix $(( v_{i,j}))$ given by :\\
$
v_{i,j}=\left\{
\begin{array}{ll}
t_i & \mbox {if $i=j$}\\
\min(t_i,t_j) & \mbox{otherwise}
\end{array}
\right.
$
In fact, the above finite dimensional convergence follows form Theorem \ref{thm.1} since the self normalised sums is converging in probability to the Wiener motion properly scaled in the sup norm metric.
\end{rmk}
\begin{rmk}
Note that we have shown finite dimensional convergence results for \\ $(S_{[nt_1]}/V_{n,p}, S_{[nt_2]}/V_{n,p},\ldots, S_{[nt_k]}/V_{n,p})$. However this is equivalent to show finite dimensional convergence for the process $Y_{n,p}(\cdot)$. To see this, note that 
\beo
E(|Y_{n,p}(t_1)- S_{[nt_1]}/V_{n,p}|^2) &=& E((nt_1-[nt_1]^2 |X_{[nt_1]}|^2/V_{n,p}^2))\\
&\le& E(|X_{[nt_1]}^2/V_{n,p}^2|) \\
&\le& E(|X_{[nt_1]}^2/V_{n,2}^2|) \ \ \forall p \le 2 \\ 
&=& \frac{1}{n} \ \mbox{since $[nt_1]<n$}
\eeo
Therefore $Y_{n,p}(t_1)-S_{[nt_1]/V_{n,p}} \stackrel{P}{\rightarrow} 0$ which implies that these two are asymptotically negligible and all limiting properties of 
$S_{[nt]}/V_{n,p}$ will be shared by $Y_{n,p}(t)$. Although we have shown the result for one dimension the result can be extended in a natural way to $k$ dimensions, ie, we can show that $(Y_{n,p}(t_1)-S_{[nt_1]/V_{n,p}}, Y_{n,p}(t_2)-S_{[nt_2]/V_{n,p}}, \ldots, Y_{n,p}(t_k)-S_{[nt_k]/V_{n,p}}) \stackrel{P}{\rightarrow} 0 $. 
\end{rmk}

\section{Tightness} 
\label{tight}
\setcounter{equation}{0}

\begin{thm}
\label{tight.thm}
The process $\{Y_{n,p}(\cdot)\}$ is tight iff $p\le \alpha \le 2$.
\end{thm}
We first prove the \textit{if} part and then the \textit{only if} part. \\
{\bf \textit{If } part:} 
The process $Y_{n,p}(\cdot)$ is tight if $p\le\alpha\le2$.

\textbf{Proof:} We first take the case that $p \le \alpha <2$. From Theorem 7.3 from \cite{Bill} the process
 $Y_{n,p}(\cdot)$ is tight iff $Y_{n,p}(0)$ is tight and for all \ $\epsilon > 0$,and $\eta>0$,\  
 $\exists \delta , (0< \delta <1)$ such that  \ $\lim_{n \rightarrow \infty} P(\omega_{Y_{n,p}}(\delta) \ge \epsilon) = 0$
  where $\omega_X(\delta)= \sup_{t-s< \delta}|X(t)-X(s)|$  is the modulas of continuity for any process
   $X(\cdot)$. Also from Equation 7.11 of \cite{Bill}  $P(\omega_X(\delta) \ge 3 \epsilon) \le \sum_{i=1}^{v} P( \sup_{t_{i-1}<s<t_i}|X(s)-X(t_{i-1})|\ge \epsilon) \;\;$ 
   for any arbitrary probability $P$ ,\ $\epsilon > 0 \ \delta > 0$ process $X(\cdot)$, \ \ and for a
  partition $0=t_0< t_1 < t_2 < \ldots < t_v = 1$ such that $\min_{1< i < v} (t_i - t_{i-1}) \ge \delta.$  \\ 
Take partition $t_i=m_i/n$ where $0=m_0 < m_1 < \ldots < m_v =n$. By the definition of the process in 
(\ref{defn eq}) we have that $\sup_{t_{i-1}<s<t_i} |Y_{n,p}(s)- Y_{n,p}(t)|=\max_{m_{i-1}<k<m_i} \frac{|S_k-S_{m_{i-1}|}}{V_{n,p}}$. 
Therefore, 
\beo
P(\omega(Y_{n,p},\delta) \ge 3 \epsilon) &\le& \sum_{i=1}^{v} P[\max_{m_{i-1}<k<m_i} |S_k - S_{m_{i-1}}|  \ge \epsilon V_{n,p}]. \\
\eeo
{The sequence $\{S_n\}$ is stationary and hence the above is same as}
\beo
& & \sum_{i=1}^{v} P[\max_{k < m_i - m_{i-1}} |S_k| > \epsilon V_{n,p}].
\eeo
Choose $m_i=mi$  where $m$ is an integer satisfying $m=\left\lceil n \delta   \right\rceil$ and $v = \left\lceil n/m  \right\rceil$. With this choice $v \rightarrow 1/\delta < 2/\delta$. Therefore for sufficiently large $n$,
\beo
P(\omega(Y_{n,p},\delta) \ge 3 \epsilon) & \le & v P(\max_{k\le m} |S_k|/V_{n,p}> \epsilon)\\
&\le& 2/\delta \;\; P(\max_{k\le m} |S_k|/V_{n,p}> \epsilon).
\eeo
Note that $S_k/V_{k,p}$ is a martingale (increments of independent mean zero random variables) and hence $|S_k|/V_{k,p}$ is a non negative sub martingale. The ratio $V_{m,p}/V_{n,p}$ has a probability limit to $(m/n)^\frac{1}{p} \rightarrow \delta^{\frac{1}{p}}$. Therefore
\beo
\frac{1}{\delta}P(\max_{k \le m} |S_k|/V_{n,p}> \epsilon) &=& \frac{1}{\delta}P(\max_{k \le m} \frac{|S_k|}{V_{k,p}} \frac{V_{k,p}}{V_{n,p}} > \epsilon)\\
&\le& \frac{1}{\delta}P(\max_{k \le m}\frac{|S_k|}{V_{k,p}}\frac{V_{m,p}}{V_{n,p}}> \epsilon).\\
\eeo
Writing $X_m=\max_{k \le m}\frac{|S_k|}{V_{k,p}}$ and $Y_m=\frac{V_{m,p}}{V_{n,p}}$ we have 
\be
\frac{1}{\delta}P(\max_{k<m}|S_k|/V_{n,p}> \epsilon)&\le& \frac{1}{\delta} P(X_m Y_m > \epsilon) \n \\
&=& \frac{1}{\delta} \{ P(X_m Y_m > \epsilon , Y_m > 2 \delta^{\frac{1}{\delta}}) + P(X_m Y_m > \epsilon , Y_m < 2 \delta^{\frac{1}{\delta}})\}\n \\
&\le& \frac{1}{\delta} \{ P(X_m Y_m > \epsilon, Y_m < 2\delta^{\frac{1}{p}}) + P(Y_m > 2 \delta^{\frac{1}{p}})\}\n \\
&\le& \frac{1}{\delta} \{ P(X_m > \epsilon/2 \delta^{\frac{1}{p}}) + P (Y_m > 2 \delta^{\frac{1}{p}})\} \n \\
&\le& \frac{1}{\delta} \{ P(X_m > \epsilon/2 \delta^{\frac{1}{p}}) + \eta\} \ \ \mbox{choosing sufficiently large m} \n \\
& & \mbox {such that $P(Y_m > 2 \delta^{\frac{1}{p}})< \eta$)}\n \\
&\le& \frac{1}{\delta} \{ 4 \delta^{\frac{2}{p}}/ \epsilon^2  V (S_m/V_{m,p}) + \eta \} \ \mbox{ \ by Doob's inequality} \n \\
& & \mbox{for nonnegative submartingales)}\n \\
&=& (4\delta^{\gamma}/\epsilon^2) V(S_m/V_{m,p}) + \eta/\delta  \;\; , \ \ \mbox{for some}\  \gamma> 0  \label{tight.2}.
\ee
Now, for $p\le \alpha <2$, or, $p < \alpha = 2, \ \ Var(S_m/V_{m,p})$ tends to zero (see section 5.1, 5.2). Tending $m \rightarrow \infty, (\mbox{since} \ \ m = \left\lceil n \delta \right\rceil)$ we have that the right hand side in (\ref{tight.2}) can be made arbitrarily small. Hence the lemma is proven. \\
For the case $p=\alpha=2$, the lemma holds by \cite{Gine} since it has been shown that the self normalised sums converges to the Normal distribution for $p=\alpha=2$.

Before proving the only if part we need the following lemma.

\begin{lem}
\label{lem.tight}
$  \{Y_{n,p}(\cdot)\}$ is tight $\Rightarrow \max_{1\le i \le n} \frac{|X_i|}{V_{n,p}} \stackrel{P}{\rightarrow} 0$.
\end{lem}

\textbf{Proof.}

We use an equivalent condition of tightness given in Theorem 4.2 of \cite{Bill}. A process is tight
 iff $\forall \epsilon >0, \forall \eta > 0, \exists n_0 $ and $0<\delta<1$ such that \\
\be
\label{tightness}
P\Bigl( \sup_{|t-s|< \delta} |Y_{n,p}(s)-Y_{n,p}(t)| \ge \epsilon \Bigr) &\le& \eta \;\; \forall t \in [0,1].
\ee	 

Assume that the hypothesis is true. Which means that for every $\epsilon, \eta >0, \exists 0<\delta<1$ such that (\ref{tightness}) holds. Choose $n_0$ sufficiently large so that $\frac{1}{n}< \delta  \;\; \forall n > n_0$. Then we have 
\beo
P(\sup_{|t-s| < \frac{1}{n}}|Y_{n,p}(t)- Y_{n,p}(s)|> \epsilon )&<& P(\sup_{|t-s| < \delta} |Y_{n,p}(t)- Y_{n,p}(s)|> \epsilon ).
\eeo
Now by definition of the process $Y_{n,p}(\cdot),$
\beo
  \sup_{|t-s|<\frac{1}{n}} |Y_{n,p}(t)-Y_{n,p}(s)| &<& \max_{1\le i \le n} \frac{|X_i|}{V_{n,p}}   \forall t \in [0,1] \\
\Rightarrow P(\max_{1\le i \le n}\frac{|X_i|}{V_{n,p}}>\epsilon) &<& P(\sup_{|t-s| < \delta} |Y_{n,p}(t)- Y_{n,p}(s)|> \epsilon )\\
\Rightarrow P(\max_{1\le i \le n}\frac{|X_i|}{V_{n,p}}>\epsilon) &<& \eta   \;\; \forall n > n_0,  \mbox{by hypothesis}. 
\eeo
\qed

\begin{rmk}
The converse is not necessarily true. To see this assume that \ 
  $\max_{1 \le i \le n}\frac{|X_i|}{V_{n,p}}\stackrel{P}\rightarrow 0$. Assume that there exists
  a $\delta_1$ such that (\ref{tightness}) holds. Given such a $\delta_1>0$, for any integer $m$ we can 
  get an $n$ such that $\frac{m}{n}< \delta_1$. Then for such a $m,n$ we have 
  $|Y_{n,p}(t)-Y_{n,p}(s)|\le (\max_{1 \le i \le n} \sum_{j=1}^{m}) |X_{i+j}|/(V_{n,p})$. But the hypothesis
  does not guarantee that the right hand side converges to zero in probability.
\end{rmk}

We use the above lemma to prove the necessary part in the following lemma.\\

{\bf \textit{Only if} part}
For $2 \ge p > \alpha$ the process is not tight.

\textbf{Proof:}\\
For $2 \ge p > \alpha$ observe that 
\beo
\max_{1\le i \le n}\frac{|X_i|}{V_{n,p}} \stackrel{P}{\rightarrow} 0\\
\Leftrightarrow \Bigl(\max_{1\le i \le n}\frac{|X_i|}{V_{n,p}}\Bigr)^p \stackrel{P}{\rightarrow} 0\\
\Leftrightarrow \max_{1\le i \le n} \frac{|X_i|^p}{\sum |X_i|^p} \stackrel{P}{\rightarrow} 0.
\eeo
But $|X_i|^p \sim DA(\gamma)$, where $\gamma=\frac{\alpha}{p} <1$, for which (Darling, \cite{Darling}, 
Theorem 5.1) says that if $Y_i \sim DA(\gamma)$ where $\gamma<1$ 
then $\max_{1\le 1 \le n}\frac{|Y_i|}{\sum |Y_i|}$ converges in distribution to a non-degenerate 
random variable $G$ whose characteristic function is identified in the same paper. 
Thus, $\max_{1\le i \le n} \frac{|X_i|^p}{\sum |X_i|^p}$  does not go to zero in probability.
Hence,
$max_{1\le i \le n}\frac{X_i}{V_{n,p}}$ cannot converge to zero in probability and therefore from Lemma \ref{lem.tight} the process cannot be tight.
\qed \\

\medskip

{
\section{Conclusion}
\label{conclusion}
The study of self normalised sums has seen a recent upsurge following the works of \cite{Gine}, \cite{Gotze}, 
\cite{Logan} and \cite{Shao}. Results for functional convergence was shown only by \cite{Csorgo}
where 
the random variable were from the domain of attraction of a Stable$(\alpha)$ 
 distribution. \\
This paper deals with the same type of random variables but with norming index $p \in {(0, 2]} $. 
Although it
 is almost intuitive that the norming index $p$ has something to do with the stability index $\alpha$ the
  relation between them has not been explored in the past. Cs\"org\"o et al \cite{Csorgo} kept the value of the 
  norming index $p$ fixed at 2 and compared with various choices of $\alpha$. This paper, to our knowledge, 
  seems to be the first one where we simultaneously change $p$ and $\alpha$. Here, using simple tools of 
  tightness and finite dimensional convergence, we show that 
  the only non trivial case is iff $p=\alpha=2$. The \textit{if} part was shown by Gine et al \cite {Gine}
  and Cs\"org\"o et al \cite{Csorgo}.
  Our paper shows the \textit{only if} part. \\
To proceed further a rate of convergence would be important. A non uniform Berry Essen bound was given in 
\cite{Bentkus}, when the random varaibles are from $DAN$, and a bound using Saddlepoint
  approximation was proved in \cite{Jing}. Although the process convergence is for $p=\alpha=2$, Logan etal
   {\cite{Logan}} have shown that the self normalised sequence can converge for $p>\alpha$. Using their
    techniques we have has shown in Section \ref{pgealpha} what the possible limiting chractetristic 
    distribution would look like. From our personal communication with
  Qi-Man Shao we have learned about an unpublished result on limiting finite-dimensional distribution
 of ($(S_{[nt_{1}]}/V_{n,p}, \ldots, S_{[nt_{k}]}/V_{n,p})$, \ $p > \alpha$) where they have shown the 
 limiting joint distribution as mixture of Poisson-type distribution using technique of 
 (Cs\"org\"o and Horvath \cite{Horvath}). The rate of convergence for this case has not been explored to our
  knowledge. 
}


\begin{thebibliography}	{10}
\bibitem {Basak}
	 Basak G.K. and Dasgupta A., (2010) {An Ornstein Uhlenbeck process associated with a Self normalising sum}.
	 Unpublished manuscript.
	 
	 \bibitem{Bentkus}
Bentkus V. and Gotze F.,(1996). The Berry-Esseen bound for Student's statistic, \textit{ Annals of Probability}, \textbf{24}, 466-490.

\bibitem{Bill}

Billingsley P. (1999) \textit{Convergence of Probability Measures}, $2^{nd}$ edition, (original edition 1968) Wiley, New York.

\bibitem{Csorgo}
	 Cs\"org\"o M., Szyszkowicz B. and Qiying W.(2003), {Donsker's theorem for Self-Normalised Partial Sum Processes}, \textit{The Annals of Probability}, 31, 3, 1228-1240
	 
	 \bibitem{Horvath}
	 Cs\"org\"o M., Horvath L. (1988), {Asymptotic representation of Self Normalized sums}, \textit{Probability and Mathematical Statistics}, 9.1, 15-24.
	 
\bibitem {Gotze}
Chistyakov G. P. and Gotze F, (2004), {Limit distributionof Studentised means}, \textit {The Annals of Probability}, \textbf{32, 1A}, 28-77. 

	 
	 \bibitem{Darling} Darling, R.(1952)
	 The influence of maximum terms in the sum of independent random variables.\textit{Transaction of the American Mathematical Society} \textbf{73, 1}, 95-107

\bibitem{erdos}
Erdos P, and Kac M., (1946), {On certain limit theorems in the theory of Probability}, \textit{Bull. of American Math. Soc.}, \textbf{52}, 292-302.

\bibitem {Feller} Feller, W. (1966) \textit{An Introduction to Probability Theory and Its Applications}, Vol. II,
 Wiley, New York. 

	 
\bibitem {Gine}
   Gine E., Gotze F. and Mason D. M. (1997), {When is a Student's t asymptotically normal}, \textit{The Annals of Probability}, \textbf{25, 3}, 1514-1531.
   


\bibitem{Jing}
Jing B. Y., Shao Qi-Man, Wang Z. (2004), Saddlepoint approximation for Students t statistics with no moment condition. \textit{The Annals of Statistics}, \textbf{32,6}, 2679-2711.

\bibitem{Kuelbs}
Griffin P. S. and Kuelbs J. D. (1989), {Self-Normalized Law of the Iterated Logarithm}, \textit{The Annals of Probability,} \textbf{17, 4}, 1571-1601.
	
\bibitem{Logan}
	 Logan B. F., Mallows . L., Rice S. O. and Shepp L. A. (1973), 
	 {Limit distribution of self normalised sums }, \textit{The Annals of Probability }, \textbf{1, 5,} 788-809.
	 
\bibitem{Obrien} Obrien G. L.(1980), {A limit theorem for sample maxima and heavy branches in Galton-Watson trees }, \textit{J. Applied Probability,} \textbf{17}, 539-545.
	 
\bibitem {Shao}
Shao, Qi-Man (1997), {Self-Normalized Large Deviations}, \textit{The Annals of Probability,} \textbf{25, 1,} 285-328.

\end{thebibliography}
\end{document}